\documentstyle[12pt,twoside]{article}
\newtheorem{theorem}{\bf Theorem}[section]
\newtheorem{proposition}[theorem]{\bf Proposition}

\newtheorem{corollary}[theorem]{\bf Corollary}

\input{amssym}

\date{}
\begin{document}

\title{{\Large\bf Weak amenability of weighted measure algebras and their second duals}}

\author{{\normalsize\sc M. J. Mehdipour and A. Rejali\footnote{Corresponding author}}}
\maketitle

{\footnotesize  {\bf Abstract.} In this paper,  we study the weak amenability of weighted measure algebras  and prove that $M(G, \omega)$ is weakly amenable if and only if $G$ is discrete and every bounded quasi-additive function is inner. We also study the weak amenability of $L^1(G, \omega)^{**}$ and $M(G, \omega)^{**}$ and show that the weak amenability of theses Banach algebras are equivalent to finiteness of $G$. This gives an answer to the question concerning weak amenability of $L^1(G, \omega)^{**}$ and $M(G, \omega)^{**}$.}
{}

\section{\normalsize\bf Introduction}

Let $G$ be a locally compact group with an identity element $e$. Let us recall that a continuous function $\omega: G\rightarrow [1, \infty)$ is called a \emph{weight} \emph{function} if for every $x, y\in G$
$$
\omega(xy)\leq\omega(x)\;\omega(y)\quad\hbox{and}\quad\omega(e)=1.
$$
Let $C_0(G, 1/\omega)$ be the
space of all functions $f$ on $G$ such that $f/\omega\in C_0(G)$, the space of all bounded continuous
functions on $G$ that vanish at
infinity. 
Let also $M(G, \omega)$ be the Banach space of all complex regular Borel measures
$\mu$ on $G$ for which $\omega\mu\in M(G)$, the measure algebra of $G$. It is well-known that
$M(G, \omega)$ is the dual
space of $C_0(G, 1/\omega)$ \cite{dl, r0, sto}, see \cite{r111, rv1} for study of weighted semigroup measure algebras; see also \cite{mr1, mr2, mr}. Note that $M(G, \omega)$ is a Banach algebra with the norm $\|\mu\|_\omega:=\|\omega\mu\|$ and the convolution product ``$\ast$"defined by 
$$
\mu\ast\nu(f)=\int_G\int_Gf(xy) \;d\mu(x) d\nu(y)\quad\quad(\mu, \nu\in M(G,\omega), f\in C_0(G, 1/\omega)).
$$
Let $L^1(G, \omega)$ be the Banach space of all Borel measurable functions $f$ on $G$ such that $\omega f\in L^1(G)$, the group algebra of $G$. Then $L^1(G, \omega)$ with the convolution product ``$\ast$" and the norm $\|f\|_{1, \omega}=\|\omega f\|_1$ is a Banach algebras.

A Borel measurable function $p$ from $G\times G$ into ${\Bbb C}$ is called \emph{quasi-additive} if for almost every where $x, y, z\in G$
$$p(xy, z)=p(x, yz)+p(y, zx).$$
If there exists $h\in L^\infty(G, 1/\omega)$ such that $$p(x, y)= h(xy)-h(yx)$$ for almost every where $x, y\in G$, then $p$ is called \emph{inner}.
Let  $D(G, \omega)$ be the set of all quasi-additive functions $p$ on $G$  such that
$$\sup_{x, y\in G}\frac{|p(x, y)|}{\omega^\otimes (x, y)}<\infty.$$
We denote by $I(G, \omega)$  the set of inner quasi-additive functions. For $\mu\in M(G, \omega)$, let $L^\infty(|\mu|, \omega)$ be the Banach space of all $\omega-$bounded Borel measurable functions $p$ on $G$ such that
$
\|p\|_{\omega, \mu}=\|p/\omega\|_\mu<\infty.
$
An element $$P=(p_\mu)_\mu\in\Pi \{L^\infty(|\mu|, \omega): \mu\in M(G, \omega)\}$$ is called a $\omega-$\emph{generlized function} on $G$ if
$$\sup\{\|p_\mu\|_{\omega, \mu}: \mu\in M(G, \omega)\}<\infty$$
and for every $\mu, \nu\in M(G, \omega)$ with $|\mu|\ll |\nu|$ we have $p_\mu=p_\nu$,$\;|\mu|-a.e.$. The space of all $\omega-$generlized function on $G$ is denoted by $GL(G, 1/\omega)$. It is well-known from \cite{sre} that $GL(G, 1/\omega)$ is the dual of $M(G, \omega)$ for the pairing
$$
\langle (p_\mu)_{\mu}, \nu\rangle=\int_Gp_\nu\;d\nu.
$$
A function $F=(F_{\mu\otimes\nu})_{\mu, \nu\in M(G, \omega)}\in GL(G\times G, \omega^\otimes)$ is called a \emph{generalized quasi-additive function} if
$$
F_{(\mu\ast\nu)\otimes\eta}(xy, z)=F_{\mu\otimes(\nu\ast\eta)}(x, yz)+F_{\nu\otimes(\eta\ast\mu)}(y, zx)
$$
for all $\mu, \nu, \eta\in M(G, \omega)$ and $x, y, z\in G$.
The set of all generalized quasi-additive functions is denoted by $GD(G, \omega)$. If there exists $p=(p_{\mu})\in GL(G, 1/\omega)$ such that for every $\mu, \nu\in M(G, \omega)$ and for almost every where $x, y\in G$
$$
F_{\mu\otimes \nu}(x, y)=p_{\mu\ast\nu}(xy)-p_{\nu\ast\mu}(yx),
$$
then $F$ is said to be a \emph{generalized inner quasi-additive function}. The set of all  generalized inner quasi-additive functions is denoted by $GI(G, \omega)$.

Let $A$ be a Banach algebra and $D: A\rightarrow A^*$ be a bounded linear operator. Then $D$ is called \emph{cyclic} if $\langle D(a), a\rangle=0$ for all $a\in A$. Let us recall that a bounded linear operator $D: A\rightarrow A^*$ is called a \emph{derivation} if $$D(ab)=D(a)\cdot b + a\cdot D(b)$$ for all $a, b\in A$. The space of all bounded continuous derivations from $A$ into $A^*$ is denoted by $\cal{Z}(A, A^*)$. If every element of $\cal{Z}(A, A^*)$ is cyclic, then $A$ is called \emph{cyclically} \emph{weakly amenable}, however, $A$ is called  \emph{weakly amenable} if every derivation $D\in\cal{Z}(A, A^*)$ is  inner; that is, there exists $z\in A^*$  such that for every $a\in A$
$$
D(a)=\hbox{ad}_z(a):=z\cdot a-a\cdot z.
$$
 The weak amenability of group algebras have been study by several authors. For example, Brown and Moran \cite{bm} studied the weak amenability of measure algebra of locally compact Abelian groups and showed that if zero is the only continuous point derivation of $M(G)$, then $G$ is discrete. Note that if $G$ is discrete, then $M(G)$ is weakly amenable, because in this case $M(G)=\ell^1(G)$ is always weakly amenable \cite{j}. One can prove that if $d$ is a non-zero continuous point derivation of $ M(G)$ at $$\varphi\in\Delta(M(G))\cup\{0\},$$ then the map $\mu\mapsto d(\mu)\varphi$ is a continuous non-inner derivation from $M(G)$ into $M(G)^*$. In other words, $M(G)$ is not weakly amenable. These facts give rise to the conjecture that for a locally compact group $G$, the Banach algebra $M(G)$ is weakly amenable if and only if $G$ is discrete; or equivalently, zero is the only continuous point derivation of $M(G)$ at a character.  Dales, Ghahramani and Helemskii \cite{dgh} proved this conjecture. Some authors investigated the weak amenability of the second dual of Banach algebras. For instance, Ghahramani, Loy and Willis \cite{glw}  proved that if $G$ is a locally compact Abelian group and $L^1(G)^{**}$ is weakly amenable, then $G$ is discrete. Forrest \cite{f} investigated the weak amenability of the dual of a topological introverted subspace $X$ of $VN(G)$. Under certain conditions, he showed that if $A(G)^{**}$ is weakly amenable, then every Abelian subgroup of $G$ is finite. As a consequence of this result, he improved the result of Ghahramani, Loy and Willis. In fact, for a locally compact Abelian group $G$, he proved that weak amenability of $L^1(G)^{**}$ is equivalent to the finiteness of $G$. Lau and Loy \cite{ll} considered  a left introverted subspace of $L^\infty(G)$ containing $AP(G)$, say $X$, and studied weak amenability of $X^*$. One can obtain the result of weak amenability of $L^1(G)^{**}$  from Lau-Loy's theorem. Finally, Dales, Lau and Strauss \cite{dls} proved that $L^1(G)^{**}$ is weakly amenable if and only if there is no non-zero continuous point derivation of $L^1(G)^{**}$ at the discrete
augmentation character; or equivalently, $G$ is finite.

This paper is organized as follow.  In Section 2 we study the weak amenability of $M(G, \omega)$  and show that  $M(G, \omega)$ is weakly amenable if and only if $G$ is discrete and every bounded quasi-additive function is inner.  We also prove that cyclic weak amenability and point amenability of $M(G, \omega)$
are equivalent to weak amenability of it. Section 3 is devoted to studies of the weak amenability of second dual of $L^1(G, \omega)$ and $M(G, \omega)$. We proved that $L^1(G, \omega)^{**}$ is weakly amenable if and only if $M(G, \omega)^{**}$ is weakly amenable; or equivalently, $G$ is finite. We verify
that cyclic weak amenability and point amenability of $L^1(G, \omega)^{**}$  and $M(G, \omega)^{**}$ are
equivalent to finiteness of $G$.

\section{\normalsize\bf Weighted measure algebras}

 Let $\omega_i$ be a weight function on  locally compact group $G_i$ for $i=1, 2$. Define the weight function $ \omega_1\otimes\omega_2$ on $G_1\times G_2$ by  $$\omega_1\otimes\omega_2(x_1, x_2)=\omega_1(x_1)\omega_2(x_2)$$ for all $x_1\in G_1$ and $x_2\in G_2$.
In the case where, $G_1=G_2=G$ and $\omega_1=\omega_2=\omega$, we set $\omega^\otimes=\omega_1\otimes\omega_2$.
The following result is needed to prove our results.

\begin{proposition}\label{tap}  Let $\omega_i$ be a weight function on  locally compact group $G_i$ for $i=1, 2$. Then
$$
M(G_1, \omega_1)\hat{\otimes} M(G_2, \omega_2)=M(G_1\times G_2, \omega_1\otimes\omega_2).
$$
\end{proposition}
{\it Proof.} Let $\eta_i\in M(G_i, \omega_i)$, for $i=1, 2$. Then for every $f\in C_0(G_1\times G_2)$, we have
$$
\langle\eta_1\otimes\eta_2, f\rangle=\int_{G_1}\int_{G_2}f(x, y)\;d\eta_1(x)d\eta_2(y).
$$
It is easy to prove that
$$
\eta_1\otimes\eta_2\in C_0(G_1\times G_2, 1/\omega_1\otimes\omega_2)^*= M(G_1\times G_2, \omega_1\otimes\omega_2).
$$
Conversely, let $\eta\in M(G_1\times G_2, \omega_1\otimes\omega_2)$. In view of Theorem Lusin's theorem, there exists sequences $(f_n)$ and $(g_n)$ in the unit ball  $C_c(G_1, 1/\omega_1)$ and $C_c(G_2, 1/\omega_2)$ with compact support, respectively, such that for almost every where $x\in G_1$ and $y\in G_2$
$$
f_n(x)\rightarrow 1\quad\hbox{and}\quad g_n(y)\rightarrow 1
$$
as $n\rightarrow\infty$. We define the functionals $\eta_1$ and $\eta_2$ by
$$
\eta_1(f)=\lim_n\eta(f\otimes g_n)\quad\hbox{and}\quad \eta_2(g)=\lim_n\eta(f_n\otimes g)
$$
for all $f\in C_0(G_1, 1/\omega_1)$ and $g\in C_0(G_2, 1/\omega_2)$. Then $\eta_1\in M(G_1, \omega_1)$, $\eta_2\in M(G, \omega_2)$. In fact,
$$
|\eta_1(f)|\leq\|\eta\|\|f\|_{\infty, 1/\omega}\quad\hbox{and}\quad|\eta_2(g)|\leq\|\eta\|\|g\|_{\infty, 1/\omega}.
$$
On the other hand,
\begin{eqnarray*}
\eta_1\otimes\eta_2(f\otimes g)&=&\lim_n\eta(f\otimes g_n)\eta(f_n\otimes g)\\
&=&\lim_n\int_{G_1\times G_2} f(x)f_n(x)g(y)g_n(y)\;d\eta(x, y).
\end{eqnarray*}
Since $\eta$ is bounded, it follows from Lebesgue dominated convergence theorem that $\|f_n\otimes g_n\|_{\infty, 1/\omega}\leq 1$ and $1\in L^1(\eta)$. Furthermore, $f_n\otimes g_n(x, y)\rightarrow 1$ for every $x\in G_1, y\in G_2$. For every $f\in C_0(G_1, 1/\omega_1)$ and $g\in C_0(G_2, 1/\omega_2)$
$$
\eta_1\otimes\eta_2(f\otimes g)=\int_{G_1\times G_2}f\otimes g(x, y)\;d\eta(x, y)=\eta(f\otimes g).
$$
It follows that for every $h\in C_0(G_1, 1/\omega_1)\otimes C_0(G_2, 1/\omega_2)$
$$
\eta_1\otimes\eta_2(h)=\eta(h)
$$
and so for every $h\in C_0(G_1\times G_2, 1/\omega_1\otimes\omega_2)$
$$
\eta_1\otimes\eta_2(h)=\eta(h).
$$
Therefore, $\eta_1\otimes\eta_2=\eta$.$\hfill\square$\\

For every $f\in L^1(G, \omega)$, we define the seminorm $T_f: M(G, \omega)\rightarrow [0, \infty)$ by $$T_f(\mu)=\|f\ast \mu\|_{1, \omega}+\|\mu\ast f\|_{1, \omega}.$$ The locally convex topology defined by the family of seminorms $(T_f)_{f\in L^1(G, \omega)}$ is called the \emph{strict topology} on $M(G, \omega)$ with respect to $L^1(G, \omega)$ (or briefly strict topology).

\begin{proposition} Let $G$ be a locally compact group  and $\omega$ be a weight function on $G$. If $p\in D(G, \omega)$, then there exists a unique bounded derivation $D\in {\cal Z}( M(G, \omega), M(G, \omega)^*)$ such that $p(x, y)=\langle D(\delta_x), \delta_y\rangle$ for all $x, y\in G$, where $\delta_\cdot$ is the Diract measure at $\cdot$.
\end{proposition}
{\it Proof.} Let $p\in D(G, \omega)$. Then $\Gamma(D_1)=p$ for some $D_1\in {\cal Z}( L^1(G, \omega), L^\infty(G, 1/\omega))$. By Proposition 2.1.6 \cite{r}, there exists $D_2\in {\cal Z}( M(G, \omega), L^\infty(G, 1/\omega))$ such that $D_2$ is strict-weak$^*$ continuous and $D_2|_{L^1(G, \omega)}=D_1$. Hence for every $f\in L^1(G, \omega)$,
\begin{eqnarray*}
\langle D_2(\delta_x), f\rangle&=&\lim \langle D_1(e_\alpha\ast\delta_x), f\rangle\\
&=&\lim \int_G\int_G p(x, y) (e_\alpha\ast\delta_x)(z) f(z)\;dzdy\\
&=&\int_G\int_G p(x, y) e_\alpha(zx^{-1}) f(y)\;dzdy.
\end{eqnarray*}
On the other hand, there exists a linear functional $T_1:L^1(G\times G,\omega^\otimes )\rightarrow{\Bbb C}$ such that
$$
\langle T_1, f\otimes g\rangle=\langle D_1(f), g\rangle
$$
for all $f, g\in L^1(G, \omega)$. Since $L^1(G\times G, \omega^\otimes )$ is a closed ideal in $M(G\times G, \omega^\otimes )$, it follows that $T_1$ has a strict continuous extension, say
$T_2:M(G, \omega)\hat{\otimes}M(G, \omega)\rightarrow{\Bbb C}$. Define $D: M(G, \omega)\rightarrow M(G, \omega)^*$ by $$\langle D(\mu), \nu\rangle=\langle T_2, \mu\otimes\nu\rangle$$ for all $\mu, \nu\in M(G, \omega)$. If $(e_\alpha)$ is a bounded approximate identity of $L^1(G, \omega)$, then for every $x\in G$,
$e_\alpha\ast\delta_x\rightarrow\delta_x$ in the strict topology. So
$$
T_2(e_\alpha\ast\delta_x\otimes e_\alpha\ast\delta_y)\rightarrow T_2(\delta_x\otimes\delta_y).
$$
Therefore,
\begin{eqnarray*}
\langle D(\delta_x), \delta_y\rangle&=&\lim \langle T_2(e_\alpha\ast\delta_x\otimes e_\alpha\ast\delta_y)\\
&=&\lim\langle D_2(e_\alpha\ast\delta_x),e_\alpha\ast\delta_y\rangle=p(x, y),
\end{eqnarray*}
as claimed.$\hfill\square$\\

In the following,  let ${\cal I}_{nn}( M(G, \omega),  M(G, \omega)^*)$ be the set of all inner derivations from $ M(G, \omega)$ into $ M(G, \omega)^*$, and let  ${\cal B}( M(G, \omega),  M(G, \omega)^*)$ be the space of bounded linear operators from $ M(G, \omega)$ into $ M(G, \omega)^*$. Define the  isometric isomorphism $\Gamma$ from Banach space ${\cal B}(M(G, \omega), M(G, \omega)^*)$ onto  $(M(G, \omega)\hat{\otimes}M(G, \omega))^*$ by
$$
\langle \Gamma(T), \mu\otimes \nu\rangle=\langle T(\mu), \nu\rangle,
$$
$ M(G, \omega)\hat{\otimes}  M(G, \omega)$ is the projective tensor product of $ M(G, \omega)$; see  Proposition 13 VI in \cite{bd}.

\begin{proposition} \label{gzms} Let $G$ be a locally compact group  and $\omega$ be a weight function on $G$. Then the following statements hold.

\emph{(i)} The function $\Gamma: {\cal Z}( M(G, \omega), M(G, \omega)^*)\rightarrow GD(G, \omega)$ is an isometric isomorphism. Furthermore, $\Gamma({\cal I}_{nn}(M(G, \omega),M(G, \omega)^*))= GI(G, \omega)$.

\emph{(ii)} If $D\in {\cal Z}( M(G, \omega), M(G, \omega)^*)$, then for every $\mu\in M(G, \omega)$ there exists $F=(F_{\mu\otimes\nu})_\nu\in GD(G, \omega)$ such that $D(\mu)=(p_{\mu, \nu})_\nu$ and $p_{\mu, \nu}(y)=\int_GF_{\mu\otimes\nu}(x, y)\;d\mu(x)$ for  almost every where $y\in G$.
\end{proposition}
{\it Proof.}  Let $D\in {\cal Z}( M(G, \omega), M(G, \omega)^*)$. Then $D\in{\cal B}(M(G, \omega), M(G, \omega)^*)$. Putting $A=B= M(G, \omega)$ in the definition of $\Gamma$, we have
$$
F:=\Gamma(D)\in (M(G, \omega)\hat{\otimes} M(G, \omega))^*=GL(G\times G, 1/\omega^\otimes )
$$
and
\begin{eqnarray}\label{zzz}
\langle D(\mu), \nu\rangle&=&\langle F, \mu\otimes\nu\rangle\nonumber\\
&=&\int_G F_{\mu\otimes\nu}(x, y)\; d(\mu\otimes\nu)(x, y)\\
&=&\int_G\int_G F_{\mu\otimes\nu}(x, y)\; d\mu(x)d\nu(y)\nonumber.
\end{eqnarray}
On the other hand, if $P=(p_\mu)_{\mu\in M(G, \omega)}$, then
\begin{eqnarray*}
\langle\hbox{ad}_P(\mu), \nu\rangle&=&\langle P\cdot\mu-\mu\cdot P, \nu\rangle\\
&=&\langle P, \mu\ast\nu\rangle-\langle P, \nu\ast\mu\rangle\\
&=&\int_Gp_{\mu\ast\nu}\;d(\mu\ast\nu)-\int_Gp_{\nu\ast\mu}\; d(\nu\ast\mu)\\
&=&\int_G\int_G(p_{\mu\ast\nu}(xy)-p_{\nu\ast\mu}(yx))\;d\mu(x)d\nu(y).
\end{eqnarray*}
Now, by the argument used in the proof of Theorem 2.3 in \cite{mr12}, it can be shown that the statement (i) holds. For (ii), assume that  $D\in {\cal Z}( M(G, \omega), M(G, \omega)^*)$ and $\mu\in M(G, \omega)$. Then $$D(\mu)\in M(G, \omega)^*=GL(G, 1/\omega).$$ Thus $D(\mu)=(p_{\mu, \nu})_{\nu}$ for some $(p_{\mu, \nu})_{\nu}\in GL(G, 1/\omega)$. Hence for every $\nu\in M(G, \omega)$, we have
$$
\langle D(\mu), \nu\rangle=\int_G p_{\mu, \nu}\; d\nu.
$$
This together with (\ref{zzz}) s
shows that $$p_{\mu, \nu}(y)=\int_GF_{\mu\otimes\nu}(x, y)\;d\mu(x)$$ for almost every where $y\in G$.$\hfill\square$\\

We are now in a position to prove the main result of this section.

\begin{theorem}\label{y} Let $G$ be a locally compact group and $\omega$ be a weight function on $G$. Then the following assertions are equivalent.

\emph{(a)} $M(G, \omega)$ is weakly amenable.

\emph{(b)} For every $D\in {\cal Z}( M(G, \omega), M(G, \omega)^*)$ there exists $P=(p_\mu)_{\mu\in M(G, \omega)}$ such that
$\langle D(\mu), \nu\rangle=\int_G\int_G(p_{\mu\ast\nu}(xy)-p_{\nu\ast\mu}(yx))\;d\mu(x)d\nu(y)$ for all $\mu, \nu\in M(G, \omega)$.

\emph{(c)} Every generalized quasi-additive function is inner.

\emph{(d)} $M(G)$ is weakly amenable and $D(G, \omega)=I(G, \omega)$.

\emph{(e)} $G$ is discrete and every non-inner quasi-additive function in $L^\infty(G, 1/\omega)$ is unbounded.
\end{theorem}
{\it Proof.} The implications (a)$\Rightarrow$(b)$\Rightarrow$(c)$\Rightarrow$(a) follow from Proposition \ref{gzms}. By Theorem 1.2 in \cite{dgh} the implication (d)$\Leftrightarrow$(e) holds. Also, the implication (e)$\Rightarrow$(a) follows from Corollary 2.5 in \cite{mr12}. For (a)$\Rightarrow$(e), let $M(G, \omega)$ be weakly amenable and $\varphi$ be a character of $M(G)$. If $d$ is a continuous  point derivation at $\varphi$ on $M(G)$, then $d|_{M(G, \omega)}$ is a continuous point derivation of $M(G, \omega)$ at $\varphi|_{M(G, \omega)}$. Hence $d$ is zero on $M(G, \omega)$. Since $M(G, \omega)$ is dense in $M(G)$, we have $d=0$ on $M(G)$ which is implies $G$ is discrete. Apply Theorem 2.4 in \cite{mr12} to conclude that $D(G, \omega)=I(G, \omega)$.$\hfill\square$\\

From Theorem 4.8  in \cite{mr4} and Theorem \ref{y} and its proof, we may prove the next result.

\begin{corollary}\label{shar} Let $G$ be a locally compact group. Then the following assertions are
equivalent.

\emph{(a)} $M(G, \omega)$ is weakly amenable.

\emph{(b)} $M(G, \omega)$ is cyclically weakly amenable.

\emph{(c)} $M(G, \omega)$ is point amenable.

\emph{(d)} $G$ is discrete and every non-inner quasi-additive function in $L^\infty(G,1/\omega)$ is unbounded.
\end{corollary}

An elementary computation shows that the functions  $\omega^\prime$ and $\omega^*$ defined by
$$
\omega^\prime(x)=\omega(x^{-1})\quad\hbox{and}\quad \omega^*(x)=\omega\otimes\omega^\prime(x,x)
$$
are weight functions on $G$.
Combining Theorem \ref{y} and the result of \cite{mr12} we have the following result.

\begin{corollary} Let $\omega$ and $\omega_0$ be weight functions on a locally compact group $G$. Then the following statements hold.

\emph{(i)} If $\omega\leq m\omega_0$ for some $m>0$, $M(G, \omega_0)$ is weakly amenable and $I(G, \omega_0)=D(G, \omega_0)$ , then $M(G, \omega)$ is weakly amenable.

\emph{(ii)} If $\omega$ and $\omega_0$ are equivalent, then weak amenability of $M(G, \omega)$ is equivalent to weak amenability of $M(G, \omega_0)$.

\emph{(iii)} $M(G, \omega^\prime)$ is weakly amenable if and only if $M(G, \omega)$ is weakly amenable.

\emph{(iv)} If $M(G, \omega^*)$ is weakly amenable and $I(G, \omega^*)=D(G, \omega^*)$, then $M(G, \omega)$ is weakly amenable.

\emph{(v)} If $G$ is Abelian, then $M(G, \omega^*)$ is weakly amenable if and only if $M(\frak{D}, \omega^\otimes)$ is weakly amenable, where $\frak{D}:=\{(x, x^{-1}): x\in G\}$.
\end{corollary}

Let $\phi: G\rightarrow G$ be a group epimorphism and $\omega$ be a weight function on $G$. Then the function $\overleftarrow{\omega}: G\rightarrow [1, \infty)$
defined by $\overleftarrow{\omega}(x)=\omega(\phi(x))$ is a weight function on $G$. For every quasi-additive function $p$, let $\frak{S}(p)$ be the quasi-additive function defined by
$$
\frak{S}(p)(x, y)=p(\phi(x), \phi(y))\quad\quad(x,y\in G).
$$
Theorem \ref{y} together with Proposition 4.1 and Theorem 4.6 in \cite{mr12} proves the next result.

\begin{corollary} Let $\omega$  be weight function on locally compact group $G$. Then the following statements hold.

\emph{(i)} If $\phi: G\rightarrow G$ is a continuous group epimorphism, $M(G, \overleftarrow{\omega})$ is weakly amenable and $\frak{S}(I(G, \omega))=I(G, \overleftarrow{\omega})$, then $M(G, \omega)$ is weakly amenable.

\emph{(ii)} If $G$ is Abelian and $M(G, \tilde{\omega})$ is weakly amenable, then $M(H, \omega|_H)$ is weakly amenable, where $H$ is a subgroup of $G$.
\end{corollary}

\begin{corollary} \label {ten1} Let $\omega_i$ be a weight function on a locally compact group $G_i$, for $i=1, 2$. Then the following assertions are equivalent.

\emph{(a}) $M(G_1, \omega_1)\hat{\otimes}M(G_2, \omega_2)$ is weakly amenable.

\emph{(b}) $M(G_1, \omega_1)\hat{\otimes}M(G_2, \omega_2)$ is  cyclically weakly amenable.

\emph{(c}) $M(G_1, \omega_1)\hat{\otimes}M(G_2, \omega_2)$ is point amenable.

\emph{(d})  $M(G_i, \omega_i)$ is weakly amenable and $G_i$ is discrete, for $i=1,2$.
\end{corollary}
{\it Proof.} Let $M(G_1, \omega_1)\hat{\otimes}M(G_2, \omega_2)$ be point amenable. Since $M(G_i, \omega_i)$ is unital, for $i=1, 2$, from Proposition \ref{tap} we infer that then $M(G_1\times G_2, \omega_1\otimes\omega_2)$ is point amenable. By Theorem \ref{y}, $G_1\times G_2$ is discrete. It follows that $G_i$ is discrete, for $i=1,2$. Hence $M(G_i, \omega_i)=\ell^1(G_i, \omega_i)$ and so $$
\ell^1(G_1, \omega_1)\hat{\otimes}\ell^1(G_1, \omega_1))=M(G_1, \omega_1)\hat{\otimes}M(G_2, \omega_2)
$$
is weakly amenable. In view of  Corollary 4.8 in \cite{mr12}, $\ell^1(G_i, \omega_i)$ is weakly amenable. So (c) implies (d).

Let $M(G_i, \omega_i)$ is weakly amenable and $G_i$ is discrete, for $i=1,2$. By Corollary \ref{shar}, $M(G_i, \omega_i)$ is point amenable, for $i=1,2$. It follows from Theorem 4.1 in \cite{mr5} and Proposition \ref{tap} that
$$
M(G_1, \omega_1)\hat{\otimes}M(G_2, \omega_2)=M(G_1\times G_2, \omega_1\otimes\omega_2)
$$
 is point amenable. Again, apply Corollary \ref{shar} to conclude that $M(G_1, \omega_1)\hat{\otimes}M(G_2, \omega_2)$ is weakly amenable. That is, (d) implies (a).$\hfill\square$\\

As a consequence of Corollary \ref{ten1}, we give the next result.

\begin{corollary}  Let $\omega_i$ be a weight function on a locally compact discrete group $G_i$, for $i=1, 2$. Then the following assertions are equivalent.

\emph{(a})  $\ell^1(G_1, \omega_1)\hat{\otimes}\ell^1(G_2, \omega_2)$ is weakly amenable.

\emph{(b})  $\ell^1(G_1, \omega_1)\hat{\otimes}\ell^1(G_2, \omega_2)$ is  cyclically weakly amenable.

\emph{(c})  $\ell^1(G_1, \omega_1)\hat{\otimes}\ell^1(G_2, \omega_2)$  is point amenable.

\emph{(d})  $\ell^1(G_i, \omega_i)$ is weakly amenable and $G_i$ is discrete, for $i=1,2$.
\end{corollary}

We say that $T\in M(G, \omega)^*$ \emph{vanishes at infinity}
if for every $\varepsilon>0$, there exists a compact
subset $K$ of $G$, for which $|\langle T, \mu \rangle|<\varepsilon$, where
$\mu \in M(G, \omega)$ with $|\mu|(K)=0$ and $\|\mu\|_\omega=1$.
We denote by $M_*(G, \omega)$ the subspace of $M(G, \omega)^*$ consisting of all functionals that
vanish at infinity. In the case where, $\omega(x)=1$ for all $x\in G$, we write
$$
M_*(G, \omega):=M_*(G).
$$
The space $M_*(G, \omega)$ is a norm closed subspace of $M(G, \omega)^*$. It is proved that $M_*(G, \omega)^*$ with the first Arens product is a Banach algebra \cite{mm}.
For each $f\in L^1(G, \omega)$, we may consider $f$ as a linear functional in $M_*(G)^*$. One can prove that $L^1(G, \omega)$ is a closed ideal in $M_*(G, \omega)^*$ and  $M_*(G, \omega)^*= L^1(G, \omega)$ if and only if $G$ is discrete \cite{mm}; see \cite{m} for the case $\omega=1$.

\begin{corollary}\label{m1} Let $G$ be a locally compact group. Then the following assertions are equivalent.

\emph{(a)} $M_*(G, \omega)^*$ is weakly amenable.

\emph{(b)} $M_*(G, \omega)^*$ is  cyclically weakly amenable.

\emph{(c)} $M_*(G, \omega)^*$ is point amenable.

\emph{(d)} $G$ is discrete and every non-inner quasi-additive function in $L^\infty(G, 1/\omega)$ is unbounded.
\end{corollary}
{\it Proof.} Let $M_*(G, \omega)^*$ be point amenable.  Since $M(G, \omega)$ is a direct summand of $M_*(G, \omega)^*$, by Theorem 3.7  in \cite{mr5}, $M(G, \omega)$ is point amenable. Hence $G$ is discrete and every non-inner quasi-additive function in $L^\infty(G, 1/\omega)$ is unbounded. Thus (c) implies (d).  It is easy to see that if $G$ discrete, then
$$
M_*(G, \omega)*=\ell^1(G, \omega)=M(G, \omega).
$$
It follows that (d) implies (a).$\hfill\square$\\

Let $L^\infty (G,1/\omega)$ be the space of all Borel measurable
functions $f$ on $G$ with $f/\omega\in L^\infty (G)$,  the Lebesgue space of bounded Borel measurable functions on $G$.  Let  also $L_0^\infty(G, 1/\omega)$ denote the
subspace of $L^\infty(G, 1/\omega)$ consisting of all functions $f\in L^\infty (G,1/\omega)$ that vanish at
infinity. It is proved that  $L_0^\infty(G, 1/\omega)$ is left
introverted in $L^\infty(G, 1/\omega)$. So  $L_0^\infty(G, 1/\omega)^*$ is a Banach algebra with the first Arens product \cite{lp}; see also \cite{dl, mmn, mnr, r0}.

\begin{corollary} Let $G$ be a locally compact group. Then the following assertions are equivalent.

\emph{(a)} $L_0^\infty(G, 1/\omega)^*$ is weakly amenable.

\emph{(b)}$L_0^\infty(G, 1/\omega)^*$ is  cyclically weakly amenable.

\emph{(c)} $L_0^\infty(G, 1/\omega)^*$ is point amenable.

\emph{(d)} $G$ is discrete and every non-inner quasi-additive function in $L^\infty(G, 1/\omega)$ is unbounded.
\end{corollary}

\begin{corollary}  Let $\omega_i$ be a weight function on a locally compact group $G_i$, for $i=1, 2$. Then the following assertions are equivalent.

\emph{(a)}  $M_*(G_1, \omega_1)^*$ and $M_*(G_2, \omega_2)^*$  are weakly amenable.

\emph{(b)} $L_0^\infty(G, 1/\omega)^*$ and $L_0^\infty(G, 1/\omega)^*$  are weakly amenable.

\emph{(c)} $M_*(G_1, \omega_1)^*\hat{\otimes}M_*(G_2, \omega_2)^*$  is weakly amenable and $G_i$ is discrete, for $i=1,2$.

\emph{(d)} $L_0^\infty(G, 1/\omega)^*\hat{\otimes}L_0^\infty(G, 1/\omega)^*$  is weakly amenable and $G_i$ is discrete, for $i=1,2$.
\end{corollary}
{\it Proof.} Assume that  $M_*(G_1, \omega_1)^*$ and $M_*(G_2, \omega_2)^*$  are weakly amenable. By Corollary \ref{m1}, $G_i$ is discrete and $M(G_i, \omega_i)=M_*(G_i, \omega_i)^*$  is weakly amenable, for $i=1,2$. It follows from Corollary \ref{ten1} that
\begin{eqnarray*}\label{ttt}
M_*(G_1, \omega_1)^*\hat{\otimes}M_*(G_2, \omega_2)^*=M(G_1, \omega_1)\hat{\otimes}M(G_2, \omega_2)
\end{eqnarray*}
is weakly amenable. So (a) implies (c).

Let $M_*(G_1, \omega_1)^*\hat{\otimes}M_*(G_2, \omega_2)^*$ be weakly amenable and $G_i$ is discrete, for $i=1,2$.
This implies that $M(G_1, \omega_1)\hat{\otimes}M(G_2, \omega_2)$ is weakly amenable. Thus  $M_*(G_i, \omega_i)^*=M(G_i, \omega_i)$  is weakly amenable. Hence (c) implies (a). Similarly, (b) and (d) are equivalent.$\hfill\square$\\

Let $\hbox{LUC}(G, 1/\omega)$ be the space of all continuous function $f$ on $G$ such that $f/\omega$ is a left uniformly continuous functions on $G$; for study of this space see \cite{za}. Let $WAP(\frak{A})$ be the  space of all weakly almost periodic  functionals on Banach algebra $\frak{A}$, that is, $f\in \frak{A}^*$ such that  the map $a\mapsto af$ from $\frak{A}$ into $\frak{A}^*$ is weakly compact, where $\langle af, b\rangle= \langle f, ba\rangle$ for all $b\in \frak{A}$.

\begin{corollary} Let $\hbox{WAP}(L^1(G, \omega))^*$ or $\hbox{LUC}(G, \omega)^*$ be 0-point amenable. Then $G$ is discrete.
\end{corollary}

Let $\frak{A}$ be one of the Banach algebras $M(G, \omega)$, $M_*(G, \omega)^*$, $L_0^\infty(G, 1/\omega)^*$,   $\hbox{WAP}(L^1(G, \omega))^*$ or $\hbox{LUC}(G, \omega)^*$.

\begin{proposition}\label{nar1} Let $G$ be a locally compact group. If $\frak{A}$ is cyclically amenable, then every element of  $CD(G, \omega)$ is inner.
\end{proposition}
{\it Proof.} Let $M(G, \omega)$ be cyclically amenable. Since $L^1(G, \omega)$ is a direct summand of $M(G, \omega)$, by Theorem 3.7 in \cite{mr5}, the Banach algebra $L^1(G, \omega)$ is cyclically amenable. It follows from Theorem 5.6 in\cite{mr12}  that every element of  $CD(G, \omega)$ is inner. For the other cases, we only need to recall that
$$
\frak{A}=M(G, \omega)\oplus\frak{B}
$$
for some closed subspace $\frak{B}$ of $\frak{A}$.$\hfill\square$

\section{\normalsize\bf The second dual of Banach algebras}

The main result of this section is the following which solves an open problem posed in \cite{ll}.

\begin{theorem}\label{lm} Let $G$ be a locally compact group. Then the following assertion are equivalent.

\emph{(a)} $L^1(G, \omega)^{**}$ is weakly amenable.

\emph{(b)} $L^1(G, \omega)^{**}$ is   cyclically weakly amenable.

\emph{(c)} $L^1(G, \omega)^{**}$ is point amenable.

\emph{(d)} $G$ is finite.
\end{theorem}
{\it Proof.}   Let  $\iota: L^1(G, \omega)\rightarrow L^1(G)$ be the inclusion map. Since $L^1(G, \omega)$ is dense in $L^1(G)$, $\iota$ is a continuous homomorphism with dense range. So $\iota^{**}: L^1(G,\omega)^{**}\rightarrow L^1(G)^{**}$ is epimorphism. Hence if $L^1(G, \omega)^{**}$ is point amenable, then by Theorem 2.1 in \cite{mr5} the Banach algebra $L^1(G)^{**}$ is point amenable. It follows that every continuous point derivation of $L^1(G)^{**}$ at the discrete augmentation character is zero. From Theorem  11.17 in \cite{dls} infer that $G$ is finite. So
(c)$\Rightarrow$(d). The implications (a)$\Rightarrow$(b)$\Rightarrow$(c) follows from Theorem 4.1 in \cite{mr4}. $\hfill\square$

\begin{corollary} Let $G$ be a locally compact group. Then the following assertion are equivalent.

\emph{(a)} $M(G, \omega)^{**}$ is weakly amenable.

\emph{(b)} $M(G, \omega)^{**}$ is cyclic  amenable.

\emph{(c)} $M(G, \omega)^{**}$ is point amenable.

\emph{(d)} $G$ is finite.
\end{corollary}
{\it Proof.}  Let $M(G, \omega)^{**}$ is point  amenable.  By Proposition 5.2 in \cite{mr5}, the Banach algebra $M(G, \omega)$ is point amenable. In view of Corollary \ref{shar}, $G$ is discrete. Hence $L^1(G, \omega)^{**}$ is weakly amenable. Now, apply Theorem \ref{lm}.$\hfill\square$\\

Let us recall that if there exists a compact invariant neighborhood of $e$ in $G$, then $G$ is called an $[IN]-$\emph{group.} The following result is an improvement of Theorem 3.4 in \cite{ll}.

\begin{theorem}\label{in} Let $G$ be a connected locally compact group. If either $G_d$
is amenable or $G$ is an $[IN]-$group, then the following assertions are equivalent.

\emph{(a)} $L^1(G, \omega)^{**}$  is weakly amenable.

\emph{(b)} $ M(G, \omega)$ is weakly amenable.

\emph{(c)} $G=\{e\}$.
\end{theorem}
{\it Proof.} Let $L^1(G, \omega)^{**}$ be weakly amenable. Since
$$
L^1(G, \omega)^{**}=M(G, \omega)\oplus C_0(G, \omega)^\perp
$$
and $C_0(G, \omega)^\perp$ is an ideal in $L^1(G, \omega)^{**}$ , we have $M(G, \omega)$ is weakly amenable. So (a)$\Rightarrow$(b).
Let's show that (b)$\Rightarrow$(c). To this end, let $M(G, \omega)$ be weakly amenable. It follows from Theorem \ref{y} that $G$ discrete and $M(G)$ is weakly amenable.
If $G_d$ is amenable, then from Theorem 3.3 in \cite{ll} we infer that $G=\{e\}$. If $G$ is an $[IN]-$group, then by Theorem 3.4 in \cite{ll}, $G$ is compact. Since $G$ is also discrete, it follows that $G$ is finite. Hence $G_d$ is amenable. Thus $G=\{e\}$. So (b)$\Rightarrow$(c). The implication (c)$\Rightarrow$(a) is clear. $\hfill\square$

\vspace{2mm}

 {\footnotesize
\noindent {\bf Mohammad Javad Mehdipour}\\
Department of Mathematics,\\ Shiraz University of Technology,\\
Shiraz
71555-313, Iran\\ e-mail: mehdipour@.ac.ir\\
{\bf Ali Rejali}\\
Department of Pure Mathematics,\\ Faculty of Mathematics and Statistics,\\ University of Isfahan,\\
Isfahan
81746-73441, Iran\\ e-mail: rejali@sci.ui.ac.ir\\

\footnotesize

\begin{thebibliography}{99}

\bibitem{bd} F. F. Bonsall and J. Duncan, Complete Normed Algebras, Ergebnisse der Mathematik und ihrer Grenzgebiete, Band 80, Springer-Verlag, Berlin/ Heidelberg/New York, 1973.

\bibitem{bm} G. Brown and W. Moran, Point derivations on $M(G)$, Bull. London Math. Soc., 8 (1) (1976) 57--64.

\bibitem{dgh} H. G. Dales, F. Ghahramani and A. Y. A. Helemskii, The amenability of measure algebras, J. London Math. Soc., (2) 66 (2002) 213--226.

\bibitem{dl} H. G. Dales and A. T. Lau, The second duals of Beurling algebras, Mem. Amer. Math. Soc., 177 (836) (2005).

\bibitem{dls} H. G. Dales, A. T. Lau and D. Strauss, Banach algebras on semigroups and on their compactifications,
Mem. Amer. Math. Soc., 205 (966) (2010).

\bibitem{f} B. Forrest, Weak amenability and the second dual of the Fourier algebra, Proc. Amer. Math. Soc., 125 (8) (1997) 2373--2378.

\bibitem{glw} F. Ghahramani, R. J. Loy and G. A. Willis, Amenability and weak amenability of second conjugate Banach algebras, Proc. Amer. Math. Soc., 124 (5)  (1996) 1489--1497.

\bibitem{j} B. E. Johnson, Weak amenability of group algebras, Bull. London Math. Soc., 23 (3) (1991) 281--284.

\bibitem{ll} A. T. Lau and R. J.  Loy,  Weak amenability of Banach algebras on locally compact groups, J. Funct. Anal., 145 (1) (1997) 175--204.

\bibitem{lp} A. T. Lau and J. Pym, Concerning the second dual of the group
algebra of a  locally compact group, J. London Math. Soc., 41 (1990) 445--460.

\bibitem{mmn} S. Maghsoudi, M. J. Mehdipour and R. Nasr-Isfahani,
Compact right multipliers on a Banach algebra
related to locally compact semigroups, Semigroup Forum, 83 (2011), no. 2, 205–213.

compact semigroups, Semigroup Forum, 83 (2) (2011) 205--213.

\bibitem{mnr} S. Maghsoudi, R. Nasr-Isfahani and A. Rejali,
Strong Arens irregularity of Beurling algebras with a locally convex topology, Arch. Math., 86 (5) (2006) 437--448.


\bibitem{mr1} S. Maghsoudi and A. Rejali, Unbounded weighted Radon measures and dual of certain function spaces with strict topology, Bull. Malays. Math. Sci. Soc., 36 (1) (2013) 211--219.

\bibitem{mr2} S. Maghsoudi and A. Rejali, On the dual of certain locally convex function spaces, Bull. Iranian Math. Soc., 41 (4) (2015) 1003--1017.

\bibitem{m} D. Malekzadeh Varnosfaderani, Derivations, Multipliers and Topological Centers of Certain
Banach Algebras Related to Locally Compact Groups, Thesis (Ph.D.)–University of Manitoba, 2017.

\bibitem{mm} M. J. Mehdipour and GH. R. Moghimi, The existence of non-zero compact right multipliers and Arens regularity of weighted Banach algebras, preprint.

\bibitem{mr} M. J. Mehdipour and A. Rejali, Regularity and amenability of weighted Banach algebras and their second dual on locally compact groups,  arXiv:2112.13286v1.

\bibitem{mr12} M. J. Mehdipour and A. Rejali, Weak amenability of weighted group algebras, arXiv:2209.08346.

\bibitem{mr4} M. J. Mehdipour and A. Rejali, Different types of weak amenability for Banach algebras,  arXiv:2209.13580.

\bibitem{mr5} M. J. Mehdipour and A. Rejali, Cohomological properties of different types of weak
amenability,  arXiv:4532770.

\bibitem{r0} H. Reiter and J. D. Stegeman, Classical Harmonic Analysis and Locally Compact Groups, London Math. Society Monographs, 22, Clarendon Press, Oxford, 2000.

\bibitem{r111} A. Rejali, The analogue of weighted group algebra for semitopological semigroups, J. Sci. Islam. Repub. Iran, 6 (2) (1995) 113--120.

\bibitem{r} A. Rejali, Weighted function spaces on topological groups, Bull. Iranian Math. Soc., 22 (2) (1996) 43--63.

\bibitem{rv1} A. Rejali and H. R. Vishki, Weighted convolution measure algebras characterized by convolution algebras, J. Sci. Islam. Repub. Iran, 19 (2) (2008) 169--173.

\bibitem{sre} Y. A. Sreider, The structure of maximal ideals in rings of measures with convolution, (Russian) Mat. Sbornik N.S., 27 (69), (1950) 297--318,  English translations 1953 in :Amer. Math. Soc. Transl. , 81, 365--391.

\bibitem{sto} R. Stokke, On Beurling measure algebras, arXiv:2107.14694v1.

\bibitem{za} Z.  Zaffar Jafar Zadeh, Isomorphisms of Banach Algebras
Associated with Locally Compact Groups, Thesis (Ph.D.)–University of Manitoba, 2015.
\end{thebibliography}
\end{document}